\documentclass[12pt,eps]{article}
\usepackage{hadronic}

\usepackage{graphicx}
\usepackage{epstopdf}
\newcommand{\ed}{\end{document}}
\newcommand{\be}{\begin{equation}}
\newcommand{\ee}{\end{equation}}
\newcommand{\bc}{\begin{center}}
\newcommand{\ec}{\end{center}}
\newcommand{\ba}{\begin{array}}
\newcommand{\ea}{\end{array}}

\newcommand{\cl}{\centerline}
\newcommand{\ts}{\textstyle}

\begin{document}
\thispagestyle{empty}

\vspace{-1cm}
\begin{center}


\vspace{1,5cm}

{\cl {\bf 
POLYNOMIAL INVARIANTS OF TORUS KNOTS AND \boldmath$(p,q)$-CALCULUS\footnote{Based on invited talks given at the 6th Petrov International Symposium on  {\it High Energy Physics, Cosmology and Gravity}, Kyiv (Ukraine), September 5 -- 8, 2013.
}  }}

\vspace{0,5cm}
{\cl{ {\bf Anatoliy M. PAVLYUK}}}
{ Bogolyubov Institute for
Theoretical Physics}\\
{ National Academy of Sciences of Ukraine}\\
{ Metrolohichna str. 14b, Kyiv 03680, Ukraine}\\
{ e-mail: pavlyuk@bitp.kiev.ua}\\
\end{center}

\vspace{1cm}
\cl{\bf  Abstract}
We introduce the deformed fermionic numbers, corresponding to the skein relations,
the main characteristics of knots and links. These fermionic numbers allow one
to restore the skein relations.
For the Alexander (Jones) skein relation we introduce  corresponding Alexander (Jones) fermionic
$q$-numbers, and for the HOMFLY skein relation -- the HOMFLY deformed $(p,q)$-numbers with
one fermionic parameter.

\vspace{0.3cm}

\vspace{0.3cm}

\newpage

\section{Introduction}

The problem of searching, investigating and physical interpreting of the skein relations for knots and links
is rather difficult one. Because of this we propose to deal with the corresponding deformed (bosonic, fermionic)
numbers, instead of working with the skein relations.

In the recent paper~\cite{Pa-13-2} we have shown that the one-parameter Alexander (Jones) skein relation is characterized
by the corresponding Alexander (Jones) bosonic $q$-numbers. The two-parameter
HOMFLY skein relation is described by the two-parameter HOMFLY  bosonic $(p,q)$-numbers. The introduced
bosonic $q$-numbers/$(p,q)$-numbers allow one to restore corresponding skein relations.

In the present paper we  show that the Alexander (Jones) skein relation can be defined
by the corresponding Alexander (Jones) fermionic $q$-numbers as well.
The HOMFLY skein relation is described by the two-parameter HOMFLY deformed $(p,q)$-numbers with one fermionic parameter.
 The introduced deformed fermionic numbers allow us to restore 
the skein relations as well.

\section{Skein Relation}

Skein relation, the most important characteristics of knots and links,  unites the
three polynomials 
$P_{L_{+}}(q),\, P_{L_{O}}(q),$  $P_{L_{-}}(q),$
and in the general form can be written as
 \be\label{skein}
P_{L_{+}}(q)=l_{1}P_{L_{O}}(q)+l_{2}P_{L_{-}}(q),\ee
where $l_{1}, l_{2}$ are coefficients.
The capital letter "L" stands for \ "Link",\ i.e. knot or link.
The $L_{+}, L_{O}, L_{-}$ denote Link with overcrossing, Link with zero crossing
and Link with  undercrossing correspondingly. 
The initial Link  $L_{+}$ turns into a simpler Link
$L_{O}$ by the surgery operation of elimination applied to 
chosen overcrossing of $L_{+}.$  The same initial Link  $L_{+}$ turns into another simpler Link
$L_{-}$ by the surgery operation of switching applied to the same
chosen overcrossing of $L_{+}$.

Applying the surgery operations to $L_{n+1, 2},$ the simplest torus knots (if $(n+1)$ is odd )  and links (if $(n+1)$ is even integer number),  
one obtains  $L_{n, 2}$ (by elimination)
and $L_{n-1, 2}$ (by switching).
Expressing it in terms of (\ref{skein}), we have the 
recurrence relation~\cite{GP1} 
\be\label{skein1}
P_{L_{n+1, 2}}(q)=l_{1}P_{L_{n, 2}}(q)+l_{2}P_{L_{n-1, 2}}(q)\,,\ee
which has the same coefficients as the skein relation~(\ref{skein}),
or in simpler form
 \be\label{skein2}
P_{n+1, 2}(q)=l_{1}P_{n, 2}(q)+l_{2}P_{n-1, 2}(q)\,.
 \ee
The axiomatic basis  of the knot theory includes the skein relation~(\ref{skein}), 
and the normalization condition
\be\label{norm}
P_{unknot}=P_{1, 2}=1.\ee

\section{Two-parameter deformed numbers}

The two-parameter $(P,Q)$-number for an integer $n$ appearing in connection with
two-parameter deformed oscillator is 
defined as~\cite{CJ} 
 \be \label{PQ}
[n]_{P,Q}={{\textstyle {P^{n}-Q^{n}}}\over{\textstyle {P-Q}}}\,. \ee
Some of the first $(P,Q)$-numbers are given here:
 \[ [1]_{P,Q}{=}1,\  [2]_{P,Q}{=}P+Q,\ [3]_{P,Q}{=}P^{2}+PQ+Q^{2},\
 [4]_{P,Q}{=}P^{3}+P^{2}Q+PQ^{2}+Q^{3}, \ldots\ .\]
The recurrence relation for $(P,Q)$-numbers
looks as  
\be\label{PQ-rec}
[n+1]_{P,Q}=(P+Q)[n]_{P,Q}-PQ[n-1]_{P,Q}\,.\ee

Comparing~(\ref{skein2}) (or directly (\ref{skein})) with~(\ref{PQ-rec}) 
allows one to introduce $q-$numbers,
corresponding to the skein relation (\ref{skein}).
For this it is necessary to find $P$ and $Q$ from
\be\label{cof}
P+Q=l_{1},\quad  PQ=-l_{2},\ee
and put them into~(\ref{PQ}).

\section{Alexander \boldmath$q$-numbers}

The Alexander skein relation~\cite{Al} 
 \be\label{alex-skein}
\Delta_{+}(q)-\Delta_{-}(q)=(q^{1\over2}-q^{-{1\over2}})\Delta_{O}(q)  
 \ee
 defines the Alexander polynomials $\Delta(q)$ 
for knots and links.
Rewriting it in the form~(\ref{skein})
 \be\label{alex-skein2}
\Delta_{+}(q)=(q^{1\over2}-q^{-{1\over2}})\Delta_{O}(q)+\Delta_{-}(q)\,
 \ee
gives the "Link coefficients" 
 \be\label{alex-l}
l_{1}^{A}=q^{1\over2}-q^{-{1\over2}},\quad  l_{2}^{A}=1\,.
 \ee
In analogy to~(\ref{skein2}),  we have the recurrence relation for Alexander polynomials of 
torus knots and links 
$L_{n, 2}:$ 
\be\label{alex-skein3}
\Delta_{n+1,2}(q)=(q^{1\over2}-q^{-{1\over2}})\Delta_{n,2}(q)+\Delta_{n-1,2}(q)\,.
 \ee
Comparing~(\ref{alex-skein3}) and~(\ref{PQ-rec}) 
one has two equations
\be\label{a1}
P+Q=q^{1\over2}-q^{-{1\over2}},\quad  PQ=-1\,,\ee
from which it follows
\be\label{a-cof}
P=q^{1\over2},\quad  Q=-q^{-{1\over2}}\,.\ee
Putting~(\ref{a-cof}) into~(\ref{PQ}),
we obtain the Alexander deformed fermionic $q-$numbers
\be \label{A-q}
[n]_{q^{{1\over2}}, -q^{-{1\over2}} }^{A}={{(q^{1\over2})^{n}-(-q^{-{1\over2}})^{n}   }
\over{q^{1\over2}+q^{-{1\over2}}    }}\,\cdot
\ee

The Alexander $q-$numbers $[n]_{q^{{1\over2}}, -q^{-{1\over2}} }^{A}$, which 
are described by the below indices, satisfy the following recurrence relation
\[  [n{+}1]_{q^{{1\over2}}, -q^{-{1\over2}} }^{A}=(q^{1\over2}-q^{-{1\over2}})[n]_{q^{{1\over2}}, -q^{-{1\over2}} }^{A}+
[n{-}1]_{q^{{1\over2}}, -q^{-{1\over2}} }^{A}, \] 
from which we find the Alexander "Link coefficients" (\ref{alex-l}), and by
putting them into (\ref{skein2}) and (\ref{skein}) we obtain (\ref{alex-skein3}) and (\ref{alex-skein2}).

The Alexander polynomial invariants for torus knots
$T(n,l)$ are given by the known formula~\cite{GP2}
\be\label{deltal}\ba{l}
 {\Delta}_{n,l}(q)=
{{(\ts q^{nl\over 2}-q^{-{nl\over
2}})(q^{1\over 2}-q^{-{1\over 2}})} \over{(\ts q^{n\over
2}{-}q^{-{n\over 2}})(q^{l\over 2}{-}q^{-{l\over 2}})}}\,,
 \ea\ee
 where $n$ and $l$ are coprime positive integers.
For $l=2\,,$ Eq. (\ref{deltal}) gives for torus knots $T(n,2):$
 \be\label{deltal2}
{\Delta}_{n,2}(q)= {{q^{n\over 2}+q^{-{n\over 2}}} \over{q^{1\over
2}+q^{-{1\over 2}}}}\,,\quad n=2m-1, 
 \ee
where $m$ are positive integers: $m=1,2,3,\ldots\,.$
Thus, from (\ref{deltal2}) and (\ref{A-q}) for torus knots $T(n,2):$
\be\label{deltal2-k}
{\Delta}_{n,2}(q)=[n]_{q^{{1\over2}}, -q^{-{1\over2}} }^{A} \,,\quad n=2m-1.
 \ee
It is easy to verify that (\ref{A-q}) in the case of even $n$ gives  the 
formula of the Alexander polynomials for torus links $L(n,2):$
\be\label{deltal2-l}
{\Delta}_{n,2}(q)=[n]_{q^{{1\over2}}, -q^{-{1\over2}} }^{A}={{q^{n\over 2}-q^{-{n\over 2}}} \over{q^{1\over
2}+q^{-{1\over 2}}}} \,,\quad n=2m.
 \ee
Thus, uniting (\ref{deltal2}) and (\ref{deltal2-l}), one has the Alexander polynomials 
for torus knots and links $L_{n,2}:$
\be\label{deltal2--l}
{\Delta}_{n,2}(q)=[n]_{q^{{1\over2}}, -q^{-{1\over2}} }^{A}
\,,\quad n=1,2,3\ldots .
 \ee

So, we introduced the Alexander fermionic $q-$numbers (\ref{A-q}), allowing to find the Alexander
skein relation (\ref{alex-skein2}).
At last, we recall that the Alexander bosonic $q-$numbers, which also lead to the Alexander
skein relation, were introduced as well~\cite{Pa-13-2}:
\be \label{A-b}
[n]_{q,q^{-1}}^{A}={{\textstyle {q^{n}-q^{-{n}}}}\over{\textstyle {q-q^{-{1}}}}}\,. \ee
They coincide with $q-$numbers of Biedenharn and Macfarlane~\cite{Bi,Ma}.
From the recurrence relation for the Alexander bosonic $q$-numbers (\ref{A-b})
\[ [n+1]_{q,q^{-1}}^{A}=(q+q^{-1})[n]_{q,q^{-1}}^{A}-[n-1]_{q,q^{-1}}^{A}  \]
one has the "knot coefficients": \[ k_{1}^{A}=q+q^{-{1}},\  k_{2}^{A}=-1\,,\]
from which the "Link coefficients" (\ref{alex-l}), defining the Alexander skein relation (\ref{alex-skein2}),
 can be found with the help of the formulas
\[ l_{2}^{A}=+(-k_{2}^{A})^{1\over2},\quad   l_{1}^{A}= +(k_{1}^{A}-2l_{2}^{A})^{1\over2}. \]

\section{Jones  \boldmath$q$-numbers} 

The Jones skein relation~\cite{Jo}   
\be\label{jones-skein}
q^{-1}V_{+}(q)-qV_{-}(q)=(q^{1\over2}-q^{-{1\over2}})V_{O}(q)
 \ee
introduces the Jones polynomials $V(q)$ for knots and links.
From (\ref{jones-skein}) 
 in the form~(\ref{skein})
 \be\label{jones-skein2}
V_{+}(q)=q(q^{1\over2}-q^{-{1\over2}})V_{O}(q)+q^{2}V_{-}(q)
 \ee
one finds the Jones "Link coefficients"
 \be\label{jones-l}
l_{1}^{V}=q^{3\over2}-q^{{1\over2}},\quad  l_{2}^{V}=q^{2}\,.
 \ee
Comparing~(\ref{jones-skein2}) and~(\ref{PQ-rec}) 
gives
\be\label{j1}
P+Q=q^{3\over2}-q^{{1\over2}},\quad  PQ=-q^{2}\,,\ee
from which 
\be\label{j-cof}
P=q^{3\over2},\quad  Q=-q^{{1\over2}}\,.\ee
Putting~(\ref{j-cof}) into~(\ref{PQ}),
we obtain the Jones fermionic $q-$numbers
\be \label{Jpq}
[n]_{q^{{3\over2}}, -q^{{1\over2}} }^{V}={{(q^{3\over2})^{n}-(-q^{{1\over2}})^{n}}\over{q^{3\over2}+q^{{1\over2}}}},
\ee
which satisfy the following recurrence relation
\[  [n{+}1]_{q^{{3\over2}}, -q^{{1\over2}} }^{V}=(q^{3\over2}-q^{{1\over2}})[n]_{q^{{3\over2}}, -q^{{1\over2}} }^{V}
+q^{2}[n{-}1]_{q^{{3\over2}}, -q^{{1\over2}} }^{V}, \]
defining  
the Jones "Link coefficients"~(\ref{jones-l}), and, therefore,
the Jones skein relation~(\ref{jones-skein2}) follows from~(\ref{skein}).

Finally, we recall that the Jones bosonic $q-$numbers~\cite{Pa-12-agg,Pa-13-2} are
\be \label{V-b}
[n]_{q^{3},q}^{V}={{\textstyle {q^{3n}-q^{{n}}}}\over{\textstyle {q^{3}-q}}}. \ee
The recurrence relation for the Jones bosonic $q$-numbers (\ref{V-b})
\[ [n+1]_{q^{3},q}^{V}=(q^{3}+q)[n]_{q^{3},q}^{V}-q^{4}[n-1]_{q^{3},q}^{V}  \]
gives the "knot coefficients": \[ k_{1}^{V}=q^{3}+q,\  k_{2}^{V}=q^{4}\,,\]
from which the "Link coefficients"~(\ref{jones-l}), leading to~(\ref{jones-skein2}), 
follow 
\[ l_{2}^{V}=+(-k_{2}^{V})^{1\over2},\quad   l_{1}^{V}= +(k_{1}^{V}-2l_{2}^{V})^{1\over2}. \]

\section{HOMFLY  \boldmath$(p,q)$-numbers}

The HOMFLY skein relation~\cite{FY}
\be\label{homfly-skein}
 p^{-1}H_{+}(p, z)-pH_{-}(p, z)=zH_{O}(p, z)
 \ee
 defines the two-parameter HOMFLY polynomials $H(p, z)$.
Let us make a change of the variable 
\[ z=q^{1\over2}-q^{-{1\over2}}. \]
Rewriting 
(\ref{homfly-skein}) in the form~(\ref{skein}), 
 \be\label{homfly-skein2}
\ H_{+}(p, q)=p(q^{1\over2}-q^{-{1\over2}})H_{O}(p, q)+p^{2}H_{-}(, q),
 \ee
 one obtains the "Link coefficients" 
 \be\label{homfly-l}
l_{1}^{H}=p(q^{1\over2}-q^{-{1\over2}}),\quad  l_{2}^{H}=p^{2}\,.
 \ee
Then from 
\be\label{j1}
P+Q=p(q^{1\over2}-q^{-{1\over2}}),\quad  PQ=-p^{2}\,,\ee
one has 
\be\label{h-cof}
P=pq^{1\over2},\quad  Q=-pq^{-{1\over2}}\,.\ee
Thus, we obtain the HOMFLY fermionic $q-$numbers
\be \label{H-aq}
[n]_{pq^{{1\over2}}, -pq^{-{1\over2}} }^{H}
=p^{n-1}\cdot{{{(q^{1\over2})^{n}-(-q^{-{1\over2}})^{n}   }\over{q^{1\over2}+q^{-{1\over2}} }   }}
=p^{n-1}\cdot [n]_{q^{{1\over2}}, -q^{-{1\over2}} }^{A}\,.
\ee
And, at the end, we remind that
the HOMFLY bosonic $q-$numbers look as
\be \label{homfly-b}
[n]^{H}_{p^{2}q,p^{2}q^{-1} }
=p^{2(n-1)}
\cdot{{\textstyle {q^{n}-q^{-{n}}}}\over{\textstyle {q-q^{-{1}}}}}\,. \ee

\end{document}